\documentclass[10pt]{amsart}

\usepackage{tikz}
\usetikzlibrary{positioning}
\usetikzlibrary{arrows}
\usepackage{amssymb,amsmath,latexsym,graphicx}
\usepackage{charter,eucal}
\usepackage{amssymb}
\usepackage{amscd}
\usepackage[all,cmtip]{xy}
\usepackage{rotating,multicol}
\usetikzlibrary{decorations.markings}

\newtheorem{theorem}{Theorem }[section]

\newtheorem{lemma}[theorem]{Lemma}
\newtheorem{observation}[theorem]{Observation}

\newtheorem{remark}[theorem]{Remark}
\newtheorem{corollary}[theorem]{Corollary}
\newtheorem{proposition}[theorem]{Proposition}
\newtheorem{principle}[theorem]{\textsc{Principle}}

\newcommand{\bt}{\begin{theorem}}
\newcommand{\et}{\end{theorem}}
\newcommand{\bmt}{\begin{maintheorem}}
\newcommand{\emt}{\end{maintheorem}}
\newcommand{\bc}{\begin{corollary}}
\newcommand{\bl}{\begin{lemma}}
\newcommand{\ec}{\end{corollary}}
\newcommand{\el}{\end{lemma}}
\newcommand{\bo}{\begin{observation}}
\newcommand{\eo}{\end{observation}}
\newcommand{\bp}{\begin{proposition}}
\newcommand{\ep}{\end{proposition}}
\newcommand{\br}{\begin{remark}}
\newcommand{\er}{\end{remark}}
\newcommand{\bpr}{\begin{principle}}
\newcommand{\epr}{\end{principle}}

\def\Aut{\mathrm{Aut}}
\def\I{\mathop{\mathrm{I}}}

\def\PSL{\mathbf{PSL}}

\def\I{\mathbf{I}}

\def\eop{\hspace*{\fill}$\blacksquare$}

\def\id{\mathrm{id}}

\def\cl{\mathrm{cl}}

\newcommand{\mC}{\mathcal{C}}

\newcommand{\mQ}{\mathcal{Q}}

\newcommand{\mA}{\mathcal{A}}
\newcommand{\mP}{\mathcal{P}}
\newcommand{\mL}{\mathcal{L}}
\newcommand{\mE}{\mathcal{E}}

\newcommand{\mW}{\mathcal{W}}

\newcommand{\mO}{\mathcal{O}}

\newcommand{\mM}{\mathcal{M}}

\newcommand{\hP}{\mathbf{P}}

\newcommand{\mS}{\mathcal{S}}

\newcommand{\mB}{\mathcal{B}}
\newcommand{\mH}{\mathcal{H}}
\newcommand{\mU}{\mathcal{U}}

\newcommand{\hS}{\mathbf{S}}

\usetikzlibrary{matrix}
\usetikzlibrary{arrows}
\usepackage{pgf}
\usepackage{lscape}
\usepackage{listings}

\usepackage{enumitem}

\title[Covers of quadrangles]{Covers of generalized quadrangles}



\author{Joseph A. Thas and Koen Thas}

\thanks{}

\address{Ghent University, Department of Mathematics, Krijgslaan 281, S22 and S25, B-9000 Ghent, Belgium}

\email{thas.joseph@gmail.com; koen.thas@gmail.com}

\date{}

\begin{document}

\maketitle

\begin{abstract}
We solve a problem posed by Cardinali and Sastry \cite{CS} about factorization of $2$-covers of finite classical generalized quadrangles. To that end, we develop a general theory of cover factorization  for generalized quadrangles, and in particular we study the isomorphism problem for such covers and associated geometries.  As a byproduct, we obtain new results about semipartial geometries coming from $\theta$-covers, and consider related problems.
\end{abstract}

\setcounter{tocdepth}{1}
\bigskip
{\footnotesize
\tableofcontents
}


\medskip
\section{Introduction}

In a recent paper \cite{CS}, Cardinali and Sastry study various aspects of geometries related to ovoids of the classical
finite generalized quadrangle $\mW(2^n)$. In particular, starting with the natural embedding $\mW(2^n) \hookrightarrow \mQ(5,2^n)$, 
they consider the geometry $\mA$ defined by $\mQ(5,2^n) \setminus \mW(2^n)$, and the geometry $\mE$ which consists of elliptic 
ovoids in $\mW(2^n)$ as points, and rosettes of such ovoids as lines (more formal details can be found in \S \ref{Def} and \S \ref{cov}). There is a natural
projection $\pi: \mA \mapsto \mE$. They then ask if any $2$-cover $\gamma: \mA \mapsto \mE$ factors through $\pi$ and an automorphism $\widetilde{\alpha}$
of $\mA$:

\bigskip
\begin{center}
\begin{tikzpicture}[>=angle 90,scale=2.2,text height=2.0ex, text depth=0.45ex]
\node (a0) at (0,3) {$\mA$};
\node (a1) [right=of a0] {$\mE$};

\node (b0) [above=of a1] {$\mA$};

\draw[->>,font=\scriptsize,thick]

(a0) edge node[above] {$\pi$} (a1);

\draw[->,font=\scriptsize,thick]

(b0) edge node[right] {$\gamma$} (a1);

\draw[->,font=\scriptsize,blue,thick,dotted,orange]
(b0) edge node[above] {$\widetilde{\alpha}$} (a0);

\end{tikzpicture}
\end{center}


This is the starting point of the present paper, and our first main result is an answer of that question in the concrete case of Cardinali and Sastry | the answer is positive. Our setting is without finiteness
restrictions, and without assuming that the quadrangles in question are classical. We just impose enough restrictions such that the geometries $\mA$ and $\mE$ exist. (The latter geometry will have as points subtended ovoids, and as lines rosettes of such ovoids.)
More precisely, we have:

\bt
\label{MHEP}
Let $\mS$ be a thick generalized quadrangle, and let $\mS' \hookrightarrow \mS$ be a thick full subquadrangle which is a geometric hyperplane of $\mS$. 
Define $\mA$ and $\mE$ as above, and let $\pi: \mA \mapsto \mE$ be the natural projection.
Then every cover $\gamma: \mA \mapsto \mE$ factorizes as $\gamma = \pi \circ \widetilde{\alpha}$, with $\widetilde{\alpha}$ an automorphism of $\mA$, if and only if every automorphism of $\mE$ is induced by an automorphism of $\mS$. 
\et

In fact, we also obtain a decomposition property at the level of the geometry $\mE$, which is used for proving the latter theorem:

\bt
\label{MLEP}
Let $\mS$ be a thick generalized quadrangle, and let $\mS' \hookrightarrow \mS$ be a thick full subquadrangle which is a geometric hyperplane of $\mS$. 
Define $\mA$ and $\mE$ as above, and let $\pi: \mA \mapsto \mE$ be the natural projection.
Then any cover $\gamma: \mA \mapsto \mE$ factorizes as $\gamma = {\alpha} \circ \pi$, with ${\alpha}$ a unique automorphism of $\mE$. 
\et

In categorical terms, this means that if $\texttt{Co}(\mA,\mE)$ is the category with objects the covers $\kappa: \mA \mapsto \mE$ and natural morphisms, then $\pi: \mA \mapsto \mE$ is an initial object. (In fact, it easily follows that each object is a zero object.)

Once both theorems are established, we obtain the affirmative answer of the question of Cardinali and Sastry as a corollary of 
Theorem \ref{MHEP}. 

After that, we proceed our study of covers. We first construct an abstract semipartial geometry from a $\theta$-cover, which generalizes a construction of Brown \cite{Brown}. From a number of general results, we obtain a second proof of Theorem \ref{MLEP}, and we develop a second approach to answer the question Cardinali and Sastry. 

In a last section, we show that the Kantor-Knuth generalized quadrangles contain subquadrangles that are $1$-subtended, an unkown result. This in its turn implies that no new semipartial geometries can arise from Brown's construction in the case of the Kantor-Knuth quadrangles. 

In an appendix we will consider a related problem, namely the extension problem of automorphisms of subquadrangles.

\medskip
\section{Definitions}
\label{Def}

\subsection{Quadrangles}

A {\em generalized quadrangle} (GQ) of {\em order $(s,t)$}, where $s$ and $t$ are cardinal numbers,
is a point-line incidence geometry with the following properties:
\begin{itemize}
\item[(a)]
Each point is incident with $t + 1$ lines and each line is incident with $s + 1$ points;
\item[(b)]
If $(x,L)$ is a non-incident point-line pair, there is precisely one point-line pair $(y,M)$ such that
$x \I M \I y \I L$ (here, ``$\I$'' denotes the incidence relation).
\item[(c)]
Two distinct points are incident with at most one line.
\end{itemize}

If both $s$ and $t$ are at least $2$, we say that the GQ is {\em thick}; otherwise it is called {\em thin}.

We use the usual notation $x \sim y$ to indicate that the points $x$ and $y$ are collinear, and any point is collinear with itself. For any point $x$, $x^{\perp} := \{ z \vert z \sim x \}$, and for any subset $Y$ of the point set, $Y^{\perp} := \cap_{y \in Y}y^{\perp}$; we make no distinction between $y^{\perp}$ and $\{y\}^{\perp}$.
In particular, we denote ${(Y^{\perp})}^{\perp}$ by $Y^{\perp\perp}$.

Finally, if $u$ and $v$ are distinct points, then $\mathrm{cl}(u,v)$ is the point set $\{ w \vert w^{\perp} \cap \{u,v\}^{\perp\perp} \ne \emptyset \}$.   

\subsection{Subquadrangles}

If $\mS = (\mP,\mL,\I)$ is a generalized quadrangle (where $\mP$ is the point set, $\mL$ the line set, and $\I$ is incidence), then a 
{\em subgeometry} is a triple $\Gamma = (\mP',\mL',\I')$ such that 
$\mP' \subseteq \mP$, $\mL' \subseteq \mL$, and $\I'$ is the induced incidence relation. 

A {\em subquadrangle} (subGQ) is a subgeometry which is a generalized quadrangle.
A subquadrangle $\mS'$ of a generalized quadrangle $\mS$ is {\em full} if for all its lines $U$, we have the property that 
any point of $\mS$ incident with $U$ is also a point of $\mS'$. Dually, we speak of ``ideal subquadrangles.''

\subsection{Ovoids}

If $\mS$ is a generalized quadrangle, an {\em ovoid} of $\mS$ is a set of points in $\mS$ such that each line 
contains precisely one point of $\mO$.

\subsection{Geometrical hyperplanes}

A {\em geometrical hyperplane} of a generalized quadrangle $\mS$ is a subgeometry $\Gamma$ each line of which contains at least $2$ points of $\mS$,  such that every line 
of $\mS$ either contains $1$ point of $\Gamma$ (and then it is not a line of $\Gamma$), or all its points are points of $\Gamma$ (and then it {\em is} a line of $\Gamma$). 

If $\mS$ is a thick generalized quadrangle of order $(s,t)$, a geometrical hyperplane is always one of the following two:
\begin{itemize}
\item
an ovoid (seen as a subgeometry);
\item
a full subquadrangle of order $(s,t/s)$, and then $t \equiv 0 \mod{s}$ if both parameters are finite.
\end{itemize}

If $\Gamma$ is a geometrical hyperplane of order $(s,1)$ of a thick generalized quadrangle of order $(s,t)$, then $t = s$.

\medskip
\section{Covers}
\label{cov}

\subsection{Morphisms}

If $\Gamma = (\mP,\mL,\I)$ and $\Gamma' = (\mP',\mL',\I')$ are point-line incidence geometries, then a {\em morphism}
$\gamma: \Gamma \mapsto \Gamma'$ is a map from $\mP \cup \mL$ to $\mP' \cup \mL'$ which sends points to points, lines to lines, and preserves incidence. If $\gamma$ is bijective and the inverse map also preserves incidence, we say $\gamma$ is an {\em isomorphism}. If $\Gamma = \Gamma'$, an isomorphism is also called an {\em automorphism}, and the set of automorphisms of $\Gamma$ naturally forms a group under composition of maps, which is denoted by $\Aut(\Gamma)$.

If $\Gamma$ is a point-line geometry and $S \leq \Aut(\Gamma)$, then for any point $v$ of $\Gamma$, $S_{[v]}$ denotes the subgroup of $S$ each element of which fixes $v$ linewise. Such groups will be needed later on.

\subsection{$\mA$ and $\mE$}

Let $\mS$ be a thick generalized quadrangle, and let $\mS' \hookrightarrow \mS$ be a thick full subquadrangle of $\mS$. 
Let $x$ be a point of $\mS \setminus \mS'$. Then $x^{\perp} \cap \mS' =: \mO_x$ is easily seen to be an ovoid of $\mS'$, which we call a {\em subtended ovoid} (by $x$). 
Let $L$ be any line of $\mS \setminus \mS'$ which meets $\mS'$ in a point $\ell$ (which is necessarily unique w.r.t. this property). Then each point 
of $L \setminus \{\ell\}$ subtends an ovoid in $\mS'$, and the set of all these ovoids is called the {\em rosette} $R_L$ of ovoids determined by $L$. 
If $\mS'$ is also a geometrical hyperplane, there is a natural bijection between the lines incident with $x$ and the points of $x^{\perp} \cap \mS'$. 

Keeping the latter hypothesis, we define the geometry $\mE$ to have as points the subtended ovoids of $\mS'$, and as lines the rosettes $R_L$. Incidence is 
symmetrized containment. We also define a geometry $\mA$, which is just the affine quadrangle which arises when taking away the geometric hyperplane
$\mS'$ of $\mS$.   

Note that there is a natural projection
\begin{equation}
\pi: \mA \mapsto \mE:  \begin{cases}
x \mapsto \mO_x \ \ \forall\  \mbox{points}\ x \\
L \mapsto R_L \ \ \forall\ \mbox{lines}\ L.
\end{cases}
\end{equation}

By the mere definition of $\mE$, $\pi$ is surjective on both points and lines.

\subsection{Covers}

Let $\gamma: \Gamma \mapsto \Gamma'$ be a morphism between point-line incidence geometries. Suppose both geometries are not empty | so they have either at least one point or at least one line. (In fact, since $\gamma$ is a map, it suffices to ask that $\Gamma$ is not empty.)
Then $\gamma$ is a {\em cover}
if $\gamma$ is locally a bijection, that is, if $x$ is any point of $\Gamma$, $\gamma$ induces a bijection between the lines incident with $x$, and the lines incident with $\gamma(x)$, and if $L$ is a line of $\Gamma$, it induces a bijection between the points incident with $L$ and the points incident with $\gamma(L)$. Sometimes we also say that $(\Gamma,\gamma)$ is a cover (of $\Gamma'$), or even that $\Gamma$ is a cover, if the covering map is clear.

In the statement of the next lemma, by ``connected point-line geometry,'' we mean a geometry with point set $\mP$ and line set $\mL$ such that any two elements in $\mP \cup \mL$ are contained in a path of elements in $\mP \cup \mL$.

\bl
If $\Gamma'$ is connected, then a cover $\gamma: \Gamma \mapsto \Gamma'$ is necessarily surjective.
\el
{\em Proof}.\quad
Let $X$ be any element of $\Gamma'$. Since $\Gamma$ is not empty, we can take some element of $\Gamma$, say $Y$. As $\Gamma'$ is connected, we then know there is a path in $\Gamma'$ in which $\gamma(Y)$ and $X$ are elements. Now apply the definition of cover to conclude that $X$ is an image for $\gamma$.
\eop \\

If $\gamma: \Gamma \mapsto \Gamma'$ is a cover, and each fiber (of lines {\em and} points) has constant size $\theta$, we say that $\gamma$, or $(\Gamma,\gamma)$, or $\Gamma$ is a {\em $\theta$-fold cover} (of $\Gamma'$) or simply {\em $\theta$-cover} (of $\Gamma'$). We also say that $\Gamma'$ is {\em $\theta$-covered} by $\Gamma$.

\br{\rm
Note that the map $\pi: \mA \mapsto \mE$ of the previous section is a cover. We call it the {\em canonical cover}.
}
\er

\medskip
\section{Factorization of morphisms, I | lower decomposition}

In this section, we aim to prove the following theorem.

\bt[Lower decomposition]
\label{MT1}
Let $\mS$ be a thick generalized quadrangle, and let $\mS' \hookrightarrow \mS$ be a thick full subquadrangle which is a geometric hyperplane of $\mS$. 
Define $\mA$ and $\mE$ as above, and let $\pi: \mA \mapsto \mE$ be the natural projection.
Then any cover $\gamma: \mA \mapsto \mE$ factorizes as $\gamma = \alpha \circ \pi$, with $\alpha$ an automorphism of $\mE$. 
\et

We call this property ``lower decomposition'' (of $\gamma$ over $\pi$, or simply of $\gamma$).

Before starting with the proof, let us remark that the assumption of $\gamma$ being a cover is essential: in general, morphisms between generalized quadrangles need not be injective | e.g.,
map every point of a quadrangle $\Gamma$ to one and the same point $x$ of a quadrangle $\Gamma'$, and every line of $\Gamma$ to one fixed line of $\Gamma'$ incident with $x$. (Morphisms preserve incidence, but not necessarily non-incidence.)\\

We will obtain the proof of Theorem \ref{MT1} in a series of short lemmas. First assume that $\gamma$ is as in the statement of Theorem \ref{MT1}; then $\alpha$ is uniquely determined by the expression
\begin{equation}
\gamma = \alpha \circ \pi
\end{equation}
if it is well defined. Namely, if $x$ is a point or a line of $\mE$, let $y \in \pi^{-1}(x)$; then $\alpha(x) = \gamma(y)$. If $\gamma(y)$ is independent of the choice of 
$y \in \pi^{-1}(x)$, then it follows easily that $\alpha$ is indeed an automorphism of $\mE$ (see below). \\

\bo
\label{pi}
If $u$ and $v$ are different points in $\mA$ for which $u^\perp \cap v^\perp = \emptyset$ in $\mA$, then $\pi(u) = \pi(v)$ (i.e., they subtend the same ovoids at infinity). Conversely, if $\pi(u) = \pi(v)$ for different points $u$ and $v$ in $\mA$, then $u^\perp \cap v^\perp = \emptyset$ in $\mA$. 
\eo
{\em Proof}.\quad
Suppose that $\{ u,v \}^{\perp} \cap \mA$ is empty.
As $\mS'$ is a geometrical hyperplane, it follows easily that, in $\mS$, we have
\begin{equation}
 \{u,v\}^{\perp} \subseteq (u^{\perp} \cap \mS') \cap (v^{\perp} \cap \mS'),
\end{equation}
hence the first statement.\\

Now suppose that
$\pi(u) = \pi(v)$ for different points $u$ and $v$ in $\mA$. If $v \in u^{\perp}$, then triangles arise in $\mS$ since they subtend the same ovoid in $\mS'$. For, if $w \in \{ u,v\}^{\perp} \cap \mA$, then the line $vw$ of $\mS$ contains a point $r$ of the ovoid, yielding a triangle with vertices $u, w, r$, contradiction.
\eop

\bc
\label{equiv}
Define the relation $\hP$ on the points of $\mA$ as follows:
\begin{equation}
x \hP y \ \ \mbox{if}\ \ x^{\perp} \cap y^{\perp} \cap \mA = \emptyset.
\end{equation}
Then $\hP$ is an equivalence relation.
\eop
\ec

\bl
\label{prop}
If $x$ is a point of $\mA$ and $y, z$ are points in $\gamma^{-1}(\gamma(x))$, then $y^{\perp} \cap z^{\perp} = \emptyset$ in $\mA$. (In other words, if $\gamma(y) = \gamma(z)$,  then $y^{\perp} \cap z^{\perp} = \emptyset$ in $\mA$.)
\el
{\em Proof}.\quad
Let $\gamma(y) = \gamma(z)$, and suppose $u$ is some point of $\mA$ such that $u \ne y \sim u \sim z \ne u$. As $\gamma$ induces a local bijection
between $u^{\perp}$ and $\gamma(u)^{\perp}$, it must follow that $\gamma(yu) \ne \gamma(zu)$. But then $\gamma(y) \ne  \gamma(z)$, contradiction. 
\eop \\

The next lemma proves the converse of Lemma \ref{prop}.

\bl
\label{wellpoint}
Let $u$ and $v$ be different points in $\mA$. 
If $u^\perp \cap v^\perp = \emptyset$ in $\mA$, then $\gamma(u) = \gamma(v)$.
\el
{\em Proof}.\quad
Suppose $u$ and $v$ are as in the statement, and suppose $\gamma(u) \ne \gamma(v)$. We consider the possibilities.\\

(1) {\em $\gamma(u) \sim \gamma(v)$}.\quad
Let $M$ be the line incident with both $\gamma(u)$ and $\gamma(v)$; as $\gamma$ is a cover, it induces a bijection between the lines incident with $u$ and 
the lines incident with $\gamma(u)$, so there is a unique line $U \I u$ such that $\gamma(U) = M$. Similarly, there is a unique line $V \I v$ such that 
$\gamma(V) = M$. As $\gamma(U) = \gamma(V)$, each point on this line has unique inverse images on both $U$ and $V$. So there is a point $w$ incident 
with $U$ such that $\gamma(w) = \gamma(v)$. But now $v^{\perp} \cap w^{\perp} \ne \emptyset$ (as $\hP$ is an equivalence relation), contradicting Lemma \ref{prop}.\\

(2) {\em $\gamma(u) \not\sim \gamma(v)$, and there is a point $z \in \mE$ such that $\gamma(u) \sim z \sim \gamma(v)$}.\quad
First note that such a point $z$ exists.
Then there are point-line pairs $(u',U)$ and $(v',V)$ such that $u \I U \I u'$ and $v \I V \I v'$ and $\gamma(u') = \gamma(v') = z$. It follows that $\gamma(U) = 
\gamma(u)z$ and $\gamma(V) = \gamma(v)z$. Some (unique) $V' \I u'$ satisfies $\gamma(V') = \gamma(V)$, and $V'$ contains a (unique) point $v''$
for which $\gamma(v'') = \gamma(v)$. As $u^{\perp} \cap v^{\perp} = \emptyset$ in $\mA$, $v '' \ne v$, and it follows that ${v''}^{\perp} \cap v^{\perp} = \emptyset$ in $\mA$ by Lemma \ref{prop}. 
By Corollary \ref{equiv}, it then follows that $u^{\perp} \cap {v''}^{\perp} = \emptyset$, contradiction. \eop \\

By Observation \ref{pi} and Lemmas \ref{prop} and \ref{wellpoint}, the map $\alpha$ is well defined on points. The proof for lines becomes easier at this point.

\bl
\label{lem4.6}
Let $L$ be a line of $\mE$. Then $\gamma$ maps any line of $\pi^{-1}(L)$ to the same line of $\mE$. Also, if $M$ is any line of $\mE$, then $\pi$ maps any line of $\gamma^{-1}(M)$ to the same line of $\mE$.
\el
{\em Proof}.\quad
Let $U, V \in \pi^{-1}(L)$, and let $v \I V$ in $\mA$. Then $r := v^{\perp} \cap \mS'$ is a point of $\mE$ on $L$; $\pi^{-1}(r)$ is a set of points in 
$\mA$ which meets every line in $\pi^{-1}(L)$ in exactly one point of $\mA$. All these points subtend $v^{\perp} \cap \mS'$. By Observation \ref{pi} and Lemma \ref{wellpoint}
$\gamma$ maps all these points to $\gamma(u)$, where $u := \pi^{-1}(r) \cap U$. It follows easily by letting $(V,v)$ vary, that all 
lines in $\pi^{-1}(L)$ are mapped on $\gamma(U)$. 

The proof of the second part of Lemma \ref{lem4.6} is similar, but now one has to rely on Observation \ref{pi} and Lemma \ref{prop}.
\eop \\

We are now ready to obtain the proof of Theorem \ref{MT1}.\\

{\em Proof of Theorem \ref{MT1}}.\quad
We have that $\alpha$ is well defined on points and lines, and it easily follows that incidence is preserved.  \eop \\

\br[The automorphism $\overline{\alpha}$]
\label{rem4.7}
{\rm
 Note that for each line $L$ of $\mE$ the lines of $\gamma^{-1}(L)$ all contain the same point $u$ of $\mS'$. 
 (For, let $U, V \in \gamma^{-1}(L)$, and suppose $U \cap \mS' \ne V \cap \mS'$; then $\pi(U) \ne \pi(V)$, so that $\gamma = \alpha \circ \pi$ implies that $\gamma(U) \ne \gamma(V)$, contradiction.) 
 If the ovoids of the rosette $L$ all share the point $u'$, then for any rosette $M \ne L$ whose ovoids share $u'$, the lines of $\gamma^{-1}(M)$ contain $u$. (For, suppose $U \in \gamma^{-1}(L)$ and $V' \in \gamma^{-1}(M)$ such that $U \cap \mS' \ne V' \cap \mS'$; then we can find a line $W$ of $\mS$ in 
 $\{U,V'\}^{\perp}$ for which the points $W \cap U =: w$ and $W \cap V' =: w'$ are points of $\mA$. The rosette defined by $W$ contains the points (subtended ovoids) $\mO_w$ and $\mO_{w'}$ of $\mE$, and so these points are collinear in $\mE$. It follows that $\gamma(w) = \alpha(\mO_w)$ and $\gamma(w') = \alpha(\mO_{w'})$ are also collinear points in $\mE$. Since $\gamma(w)$ is a point of the rosette $L$, and $\gamma(w')$ is a point of the rosette $M \ne L$, this is clearly not possible.) 
 If $u' =: \zeta(u)$,  and we let $u$ vary, we obtain a permutation
 \begin{equation}
 \zeta: w \mapsto w' 
 \end{equation}
 of the points of $\mS'$. For distinct collinear points $v$ and $v_1$ of $\mS'$, we have that $\zeta(v) =: v'$ and $\zeta(v_1) =: v_1'$ are also collinear as otherwise $v'$ and $v_1'$ belong to a common point (subtended ovoid) $\mO'$  of $\mE$, and hence $\alpha^{-1}(\mO') = \mO$ contains $v$ and $v_1$, contradiction. Similarly, if $v'$ and $v_1'$ are distinct collinear points of $\mS'$, then $\zeta^{-1}(v')$ and $\zeta^{-1}(v_1')$ are collinear as well.  
 It follows that the map $\zeta$ defines an automorphism $\overline{\alpha}$ of $\mS'$.}
 \er

The next result yields more information about lower decomposition.

\bc 
Let {\rm ${\texttt{Co}}(\mA,\mE)$} be the category with objects the covers $\kappa: \mA \mapsto \mE$ and morphisms defined as follows: 
\begin{quote}
elements of {\rm $\texttt{hom}(\kappa': \mA \mapsto \mE,\kappa'': \mA \mapsto \mE)$} are morphisms $\vartheta: \mE \mapsto \mE$ such that the following diagram commutes:
\begin{center}
\begin{tikzpicture}[>=angle 90,scale=2.2,text height=2.0ex, text depth=0.45ex]
\node (a0) at (0,3) {$\mE$};
\node (a0') [right=of a0] {};
\node (a1) [right=of a0'] {$\mE$};

\node (b0) [above=of a1] {$\mA$};

\draw[->,font=\scriptsize,thick,>=angle 45,orange]

(a0) edge node[above] {$\vartheta$} (a1);



\draw[->>,font=\scriptsize,thick]

(b0) edge node[right] {$\kappa''$} (a1);

\draw[->>,font=\scriptsize,thick]
(b0) edge node[above] {$\kappa'$} (a0);

\end{tikzpicture}
\end{center}
\end{quote}
Then $\pi: \mA \mapsto \mE$ is an initial object. It follows that each object is a zero object.
\ec

{\em Proof.}\quad
The fact that $\pi: \mA \mapsto \mE$ is initial is a direct corollary of Theorem \ref{MT1}. The fact that it is zero (so also terminal) follows from the fact that $\alpha$ (in Theorem \ref{MT1}) is an automorphism. Now consider two objects $\gamma: \mA \mapsto \mE$ and 
$\gamma': \mA \mapsto \mE$ in $\texttt{Co}(\mA,\mE)$. Using the fact that $\pi: \mA \mapsto \mE$ is a zero object, one now easily proves that there are unique $\delta: \mE \mapsto \mE$ and
$\delta' = \delta^{-1}: \mE \mapsto \mE$ such that the following diagram commutes:

\bigskip
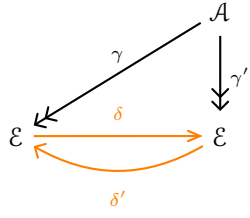
\begin{center}
\begin{tikzpicture}[>=angle 90,scale=2.2,text height=2.0ex, text depth=0.45ex]
\node (a0) at (0,3) {$\mE$};
\node (a0') [right=of a0] {};
\node (a1) [right=of a0'] {$\mE$};

\node (b0) [above=of a1] {$\mA$};

\draw[->,font=\scriptsize,thick,>=angle 45,orange]

(a0) edge node[above] {$\delta$} (a1);

\draw[->,font=\scriptsize,thick,orange]

(a1) [bend left] edge node[below] {$\delta'$} (a0);

\draw[->>,font=\scriptsize,thick]

(b0) edge node[right] {$\gamma'$} (a1);

\draw[->>,font=\scriptsize,thick]
(b0) edge node[above] {$\gamma$} (a0);

\end{tikzpicture}
\end{center}
\eop

\medskip
\bc
\label{cor2}
Let $\mS$ be a thick generalized quadrangle, and let $\mS' \hookrightarrow \mS$ be a thick full subquadrangle which is a geometric hyperplane of $\mS$. 
Define $\mA$ and $\mE$ as above, and let $\pi: \mA \mapsto \mE$ be the natural projection.
Then any $2$-cover $\gamma: \mA \mapsto \mE$ factorizes as $\gamma = \alpha \circ \pi$, with $\alpha$ some automorphism of $\mE$. 
\eop
\ec


\bc
Let $\mQ(4,q)$ be a subquadrangle of $\mQ(5,q)$. 
Define $\mA$ and $\mE$ as before, and let $\pi: \mA \mapsto \mE$ be the natural projection.
Then any $2$-cover $\gamma: \mA \mapsto \mE$ factorizes as $\gamma = \alpha \circ \pi$, with $\alpha$ some automorphism of $\mE$. 
\ec
{\em Proof}.\quad
The proof follows from Corollary \ref{cor2} and the fact that $\mQ(4,q)$ is a geometrical hyperplane of $\mQ(5,q)$.
\eop \\

\section{Factorization of morphisms, II | higher decomposition}

We now handle higher decomposition.

\subsection{Higher decomposition and higher extension}

We keep using the notation of the previous section. In that section, we have shown that given a cover $\gamma: \mA \mapsto \mE$, there exists an automorphism $\alpha$ of $\mE$ such that $\gamma = \alpha \circ \pi$, with $\pi$ the canonical cover. Now {\em suppose} that $\alpha$ itself is induced by an automorphism $\widetilde{\alpha}$ of $\mS$ (so that 
$\widetilde{\alpha} \in \Aut(\mS)_{\mS'}$, and $\widetilde{\alpha}$ extends $\overline{\alpha}$). It is then obvious that
\begin{equation}
\gamma = \alpha \circ \pi = \pi \circ \widetilde{\alpha},
\end{equation}
and so we have  (in this general setting) a ``higher decomposition'' of morphisms as asked (for every cover) by Cardinali and Sastry. Vice versa, if an automorphism $\widetilde{\beta}$ in $\Aut(\mS)_{\mS'}$ exists such that 
$\gamma = \pi \circ \widetilde{\beta}$, then $\widetilde{\beta}$ necessarily induces $\alpha$, or equivalently extends $\overline{\alpha}$, so the existence of the decomposition of Cardinali and Sastry for $\alpha$ is {\em equivalent} with the existence of an extending automorphism $\widetilde{\beta}$. In other words, we have the following result.

In the next theorem, we say that $\pi: \mA \mapsto \mE$ has the {\em higher decomposition property} if any cover $\gamma: \mA \mapsto \mE$ factorizes as $\gamma = \pi \circ \widetilde{\alpha}$ for some $\widetilde{\alpha} \in \Aut(\mS)$. By abuse of terminology, we also say that $\gamma: \mA \mapsto \mE$ has the higher decomposition property if it factorizes as above. 

\bt[Extension property]
\label{HEP}
The cover $\pi: \mA \mapsto \mE$ has the higher decomposition 
property if and only if any automorphism of $\mE$ is induced by an automorphism of $\mS$.
\eop 
\et

We call the latter property the ``higher extension property'' (for the cover $\pi: \mA \mapsto \mE$). If $\iota: \mS' \mapsto \mS$ is a natural embedding of generalized quadrangles, then we say that it has the ``higher extension property'' if each automorphism of $\mS'$ extends to an automorphism of $\mS$.
Theorem \ref{HEP} yields a criterion which makes it rather 
easy to check whether higher decomposition is true for all covers $\gamma: \mA \mapsto \mE$. We need to develop some more tools in order to apply the criterion.

\br{\rm
It is very important to note that $\widetilde{\alpha}$ need not be unique with respect to this property: {\em any} $\widetilde{\alpha}$ which extends $\overline{\alpha}$ gives rise to the same decomposition. This is in stark contrast with the automorphism $\alpha$, which is necessarily unique. If $\widetilde{\alpha}'$ is another such extending automorphism, $\widetilde{\alpha}{(\widetilde{\alpha}')}^{-1}$ fixes $\mS'$ elementwise, and it easily follows that we have a natural bijection
\begin{equation}
\eta: \Aut(\mS)_{[\mS']} \mapsto E(\overline{\alpha}): \theta \longrightarrow \widetilde{\alpha}\theta,
\end{equation}
where $E(\overline{\alpha})$ is the set of automorphisms of $\mS$ that extend $\overline{\alpha}$.}
\er

Let $\Aut(\mS)_{\mS'}$ be the stabilizer of $\mS'$ in $\mS$, and let $N = \Aut{\mS}_{[\mS']}$ be the subgroup of elements that fix $\mS'$ elementwise. Then obviously $\Aut(\mS)_{\mS'}/N \leq \Aut(\mS')$. Each element  of $\Aut(\mS)_{\mS'}$ induces an element of $\Aut(\mE)$ in a unique way, and modulo $N$ this happens in a faithful way. So 
\begin{equation}
\Aut(\mS)_{\mS'}/N \leq \Aut(\mE). 
\end{equation}

The following is immediate.

\bo
\label{obs}
If every automorphism of $\mS'$ extends to an automorphism of $\mS$, then every automorphism of $\mS'$ induces an automorphism of $\mE$ in a faithful manner.
\eop
\eo

Now if $\Aut(\mS)_{\mS'}/N \cong \Aut(\mS')$ as in Observation \ref{obs}, then the higher decomposition property is not true if and only if $\Aut(\mS') \ne \Aut(\mE)$.

\medskip
\section{The original question}

When $\mS = \mQ(5,q)$ and $\mS' = \mQ(4,q)$, the higher decomposition property is always true (in any characteristic), since any automorphism of $\mQ(4,q)$ extends in two ways (since the group of automorphisms of $\mQ(5,q)$ fixing $\mQ(4,q)$ elementwise has order $2$). In Brown \cite[Theorem 3.3]{Brown}, the following is shown:

\bt[Brown]
\label{thmbrown}
If $\Gamma$ is a finite thick GQ of order $(u,u^2)$ and $\Delta$ a subGQ of order $u$ which is $2$-covered by $\Gamma$, 
and $\Gamma'$ is a GQ of order $(u,u^2)$ with a subGQ $\Delta'$ of order $u$ which is also $2$-covered by $\Gamma'$, then the geometry $\mE$ is isomorphic to the geometry $\mE'$  if and only if there is an isomorphism $\kappa: \Delta \mapsto \Delta'$ that induces an isomorphism $\kappa': \mE \mapsto \mE'$.  
\et

As a corollary of the proof of Theorem \ref{thmbrown}, he deduces the following.

\bc[Brown]
\label{corbrown}
If $\Gamma$ is a finite thick GQ of order $(u,u^2)$ and $\Delta$ a subGQ of order $u$ which is $2$-covered by $\Gamma$, then $\Aut(\mE)$ is the stabilizer of $\mE$ in $\Aut(\mS')$. 
\ec

\br{\rm
We will generalize Theorem \ref{thmbrown}  later on. Note that Theorem \ref{thmbrown} also follows from the same reasoning as in Remark \ref{rem4.7}.}
\er

Combining Theorem \ref{HEP}, Observation \ref{obs} and Corollary \ref{corbrown} we conclude that the question of Cardinali and Sastry has an affirmative answer.


\medskip
\section{More on subtended ovoids}

In this section, as in the sections \ref{geom}--\ref{unique}--\ref{part}--\ref{last}, we only consider finite GQs.

Let $\mS' = (\mP',\mB',\I')$ be a proper subquadrangle of order $(s,t')$, $s \ne 1$, of a generalized quadrangle $\mS = (\mP,\mB,\I)$ of order $(s,t)$. Then $t \geq st'$ and each point $x$ of $\mS$ not in $\mS'$ is collinear with exactly $st' + 1$ points of an ovoid $\mO_x$ of $\mS'$; the ovoid $\mO_x$ is said to be {\em subtended} by $x$. 
The points of $\mS$ not in $\mS'$ are called the {\em external points} of $\mS'$. If $t = st'$, then each line of $\mS$ is incident with a point of $\mS'$. 

Suppose that the subtended ovoid $\mO_x$ of $\mS'$ is subtended by $\theta$ points of $\mP \setminus \mP'$, $\theta > 1$; then we say that $\mO_x$ is {\em $\theta$-subtended}. By \cite[\S 1.4.1]{PT2} we have that 
\begin{equation}
(\theta - 1)t' \leq s.
\end{equation}
If $\theta = s + 1$, $t' = 1$, and if $S_x$ is the set of the $s + 1$ points subtending $\mO_x$, then any point $z \in \mP \setminus (\mO_x \cup S_x)$ is collinear with exactly two points of $\mO_x \cup S_x$ \cite[\S 1.4.1(i)]{PT2}. If $(\theta - 1)t' = s$, $t' \ne 1$ and $\mO_x = \mO_{x'}$ with $x \ne x'$, then by \cite[\S 1.4.1(ii)]{PT2} we have that $t = st'$, and
each point $w \not\in \cl(x,x')$ is collinear with $t/s + 1 = t' + 1$ points of $\mO_x$.  \\

We now provide some examples of GQs $\mS, \mS'$ with $\mS' \subset \mS$, where each subtended ovoid $\mO_x$ is $\theta$-subtended for some $\theta$.

\subsection{Examples}
\label{ex1}

\begin{itemize}
\item[(1)]
Consider the subquadrangle $\mQ(4,q)$ of the GQ $\mQ(5,q)$. Then $s = t' = q$, $t = q^2$, and $\theta = 2$. Here $t = st'$, $(\theta - 1)t = s^2$ and $(\theta - 1)t' = s$.
\item[(2)]
Consider the subquadrangle $\mQ(3,q)$ of the GQ $\mQ(5,q)$. Then $s = q$, $t = q^2$, $t' = 1$, and $\theta = q + 1$. Here $t \ne st'$, $(\theta - 1)t \ne s^2$ and $(\theta - 1)t' = s$.
\item[(3)]
Consider the subquadrangle $\mQ(3,q)$ of the GQ $\mQ(4,q)$. Then $s =  t = q$ and $t' = 1$. If $q$ is odd,  then $\theta = 2$. For $q$ even each subtended ovoid of $\mQ(3,q)$ is subtended by just one point of $\mQ(4,q)$.
\item[(4)]
Consider the subquadrangle $\mH(3,q^2)$ of the GQ $\mH(4,q^2)$. Then $s = q^2$, $t = q^3$, $t' = q$, and $\theta = q + 1$. Here $t = st'$, $(\theta - 1)t = s^2$ and $(\theta - 1)t' = s$.
\item[(5)]
It is well known that the GQ of order $(q,q^2)$ arising from the Kantor-Knuth flock has subquadrangles of order $q$ isomorphic to $\mQ(4,q)$. This also yields examples
with $s = t' =q$, $t = q^2$, and $\theta = 2$; see \cite[\S 5]{Brown}. 
\end{itemize}

\br{\rm
Let $\mS'$ be a subquadrangle of order $s$ of the GQ $\mS$ of order $(s,s^2)$, $s \ne 1$. If every subtended ovoid of $\mS'$ is $2$-subtended, then by \cite{long} we have:
\begin{itemize}
\item[(i)]
for $s$ even we have $\mS' \cong \mQ(4,s)$, and
\item[(ii)]
for $s$ odd all points of $\mS'$ are antiregular.
\end{itemize}
}
\er

\bl
\label{7.2}
If $t' \ne 1$, $t = st'$, $(\theta - 1)t = s^2$, if every subtended ovoid of $\mS'$ is $\theta$-subtended, and if $x' \in \mP \setminus \mP'$ is collinear with no point of $\mO_x^{\perp}$, then $\vert \mO_x \cap \mO_{x'} \vert = t/s  + 1$. If $x' \not\in \mO_x^{\perp}$ is collinear with a point of $\mO_x^{\perp}$, then $\vert \mO_x \cap \mO_{x'} \vert = 1$.
\el
{\em Proof.}\quad
The second part of the statement is clear. So assume that $x' \in \mP \setminus \mP'$ is collinear with no point of $\mO_x^{\perp}$. Let $y \in \mO_x^{\perp} \setminus \{x\}$. As $x' \not\in \cl(x,y)$, $\vert \cl(x,y) \vert = \theta$ and $(\theta - 1)t = s^2$, the point $x'$ is collinear with $t/s + 1$ points of $\{ x,y \}^{\perp} = \mO_x$, and so
$\vert \mO_x \cap \mO_{x'} \vert = t/s + 1$. \eop \\

\medskip
\section{The geometry $\Gamma(\mS,\mS') = \mE$}
\label{geom}

Assume again that every subtended ovoid $\mO_x$ of $\mS'$ is $\theta$-subtended, $\theta > 1$. Let $L$ be a line of $\mS$ not in $\mS'$ which has a 
point $z$ in common with $\mS'$. If $x_1,x_2,\ldots,x_s$ are the points of $L$ not in $\mS'$, then the ovoids $\mO_{x_1},\mO_{x_2},\ldots,\mO_{x_s}$ form a 
rosette $R_L$ of ovoids in $\mS'$. Any two elements of $R_L$ have just $z$ in common and $\vert R_L \vert = s$. Also, the union of all ovoids of $R_L$
is $\mP'$. Now we introduce the point-line geometry $\mE = (\Omega,\Phi,\in)$, with $\Omega$ the set of all subtended ovoids of $\mS'$ and with $\Phi$ the set
of all rosettes of subtended ovoids in $\mS'$.

\bt
Let $\mS$ be a GQ of order $(s,t)$ and $\mS'$ a subGQ of order $(s',t')$.
If every subtended ovoid $\mO_x$ of $\mS'$ is $\theta$-subtended, $\theta > 1$, with $t' \ne 1$, $t = st'$, and $(\theta - 1)t = s^2$, then the geometry $\mE = (\Omega,\Phi,\in)$ is a semi partial geometry with parameters 
$$s^* = s - 1,\ \ t^* = t,\ \ \alpha^* = \theta, \ \ \mu^* = \theta(t - t').$$
\et
{\em Proof}.\quad
Every line of $\mE$ contains $s$ points of $\mE$, every point of $\mE$ is contained in $st' + 1 = t + 1$ lines of $\mE$, and from Lemma \ref{7.2} it easily  follows that any two distinct  points 
of $\mE$ are incident with at most one line of $\mE$. 

Ler $R_L \in \Phi$, $\mO_x \in \Omega$ and $\mO_x \not\in R_L$. Let $z$ be the common point of the ovoids in $R_L$, and let $x = x_1,x_2,\ldots,x_\theta$ be the elements of $\mS$ subtending $\mO_x$. If $x$ and $z$ are collinear in $\mS$, then also $x_i$ and $z$, with $i = 1,2,\ldots,\theta$, are collinear in $\mS$. In such a case the point $\mO_x$ of $\mE$ is collinear with no point of the line $R_L$ of $\mE$. Now assume that $x$ and $z$ are not collinear in $\mS$; then $x_i$ and $z$, with $i = 1,2,\ldots,\theta$, are not collinear in $\mS$. Let $u_i \sim x_i$, with $u_i \I L$, $i = 1,2,\ldots,\theta$. As $\{ x_i,x_j \}^{\perp} = \mO_x$, $i \ne j$, we have $u_i \ne u_j$. If $M_i = x_iu_i$, with $i = 1,2,\ldots,\theta$, then the rosettes $R_{M_i}$ are the lines of $\mE$ which contain $\mO_x$ and are concurrent with $R_L$ in $\mE$. Remark that these $\theta$ rosettes are independent of the choice of the line $L'$ with $R_L = R_{L'}$. Consequently, $\alpha^* = \theta$.

Now let $\mO_x$ and $\mO_{x'}$ be two points of $\mE$ which do not belong to a common rosette of $\Phi$. Let $x = x_1,x_2,\ldots,x_{\theta}$ be the elements of $\mS$ subtending $\mO_{x}$ and let $x' = x_1',x_2',\ldots,x_{\theta}'$ be the elements of $\mS$ subtending $\mO_{x'}$. The points of $\{ x_i,x_j'\}^{\perp}$ in $\mS$ but not in $\mS'$, with $i, j =1,2,\ldots,\theta$, subtend the ovoids collinear with $\mO_x$ and $\mO_{x'}$ in $\mE$. Clearly, we may choose $x_i$, say $x_i = x$. By Lemma \ref{7.2} we have $\vert \mO_x \cap \mO_{x_j'} \vert = \vert \mO_x \cap \mO_{x'} \vert = t/s + 1 = t' + 1$, with $j = 1,2,\ldots,\theta$. Hence 
\begin{equation}
\vert \Big(   \{ x,x_1' \}^{\perp} \cup \ldots \cup \{ x,x_{\theta}'\}^{\perp}   \Big) \setminus \mP'  \vert = \theta(t - t'). 
\end{equation}

We conclude that $\mu^* = \theta(t - t')$, and this proves the theorem.
\eop \\

\subsection{Examples}

\begin{itemize}
\item[(1)]
Considering the items (1) and (5) of \S\S \ref{ex1}, we obtain semi partial geometries with parameters
\begin{equation}
s^* = q - 1, \ \ t^* = q^2, \ \ \alpha^* = 2, \ \ \mu^* = 2q(q - 1).
\end{equation}
\item[(2)]
Considering the item (4) of \S\S \ref{ex1}, we obtain a semi partial geometry with parameters
\begin{equation}
s^* = q^2 - 1, \ \ t^* = q^3, \ \ \alpha^* = q + 1, \ \ \mu^* = q(q + 1)(q^2 - 1).
\end{equation}
\end{itemize}

\br{\rm
For the geometries described in the previous examples, we also refer to Brown \cite{Brown} and Hirschfeld and Thas \cite{HT2}.
}
\er

\medskip
\section{Uniqueness of covers}
\label{unique}

Consider again the geometry $\mE = \Gamma(\mS,\mS') = (\Omega,\Phi,\in)$. Then $\mE$ is {\em $\theta$-covered} by the geometry $\mA = 
\mS \setminus \mS' = (\mP \setminus \mP',\mB \setminus \mB',\I'')$, with $\I''$ the restriction of $\I$ to $\Big((\mP \setminus \mP') \times (\mB \setminus \mB')\Big) \cup \Big((\mB \setminus \mB') \times (\mP \setminus \mP')\Big)$; the geometry $\mA$ through the mapping $\gamma:\mA \mapsto \mE$ is a {\em $\theta$-cover} of $\mE$. 

The following theorem generalizes Theorem 3.3 of \cite{Brown}.

\bt
If $\gamma: \mC \mapsto \mE$ is any $\theta$-cover of $\mE$ without triangles, then there exists a GQ $\chi$, a subquadrangle $\chi'$ of $\chi$ and an isomorphism $\sigma^*$ from $\chi'$ to $\mS'$, such that $\chi \setminus \chi' = \mC$ and $\Gamma(\chi,\chi')^{\sigma^*} = \mE$.
\et

{\em Proof}.\quad
It is clear that each line of $\mC$ has $s$ points and that each point of $\mC$ is on $t + 1$ lines. 

Let $x$ be a point of $\mS'$. The lines of $\mS$ not in $\mS'$ containing $x$ define $(t - t')/\theta$ rosettes by canonical projection, and the set of these rosettes is covered by a set of $t - t'$ mutually disjoint lines of $\mC$. The set of these $t - t'$ lines of $\mC$ is by definition a point $x^*$ of $\chi'$. If we consider $s + 1$ collinear points $x_1,x_2,\ldots,x_{s + 1}$ of $\mS'$, then the corresponding points $x_1^*,x_2^*,\ldots,x_{s + 1}^*$ constitute a line $M^*$ of $\chi'$. In such a case the lines of the corresponding $s + 1$ sets of $t - t'$ lines of $\mC$ are mutually disjoint. Each point of $\chi'$ is contained in $t' + 1$ lines of $\chi'$. The incidence structure $\chi'$ is isomorphic to the GQ $\mS'$, and 
\begin{equation}
\sigma^*: x^* \mapsto x
\end{equation}
defines an isomorphism from $\chi'$ to $\mS'$. 

Let us now extend the geometry $\mC$ to a geometry $\chi$. The line set of $\chi$ is the union of the line set of $\mC$ and the line set of $\chi'$; the point set of $\chi$ is the union of the point set of $\mC$ and the point set of $\chi'$. A line $N$ of $\mC$ is incident with a point $x^*$ of $\chi'$ if $N$ belongs to the set $x^*$; the other incidences are the incidences of $\mC$ and $\chi'$. Each line of $\chi$ contains $s + 1$ points of $\chi$, and each point of $\chi$ is contained in $t + 1$ lines of $\chi$. The geometry $\chi$ contains $(s + 1)(st + 1)$ points and $(t + 1)(st + 1)$ lines.

Assume, by way of contradiction, that the non-collinear points $z_1, z_2, z_3$ of $\chi$ form a triangle. As we assumed that $\mC$ does not contain triangles, at least one of the points $z_1, z_2, z_3$ belongs to $\chi'$; also at least one of the points $z_1,z_2,z_3$ belongs to $\mC$ (as $\chi'$ has no triangles as well). Let $z_1, z_2$ belong to $\mC$, and $z_3$ to $\chi'$. The points $z_1, z_2$ cover ovoids $\mO_{u_1}, \mO_{u_2}$ of $\mS'$ with $\mO_{u_1} \cap \mO_{u_2} = \{ u\}$, and with $z_3$ corresponds a point $\widetilde{z_3}$ of $\mS'$ which belongs to $\mO_{u_1}$ and $\mO_{u_2}$. So $u = \widetilde{z_3}$ and the line $z_1z_2$ of $\chi$ contains the point $z_3$ of $\chi$, a contradiction. Next, let $z_1, z_2$ belong to $\chi'$, and $z_3$ to $\mC$. With $z_3$ corresponds an ovoid $\mO_{u_3}$ of $\mS'$, and with the lines $z_1z_3$ and $z_2z_3$ correspond rosettes containing $\mO_{u_3}$. Hence with $z_1, z_2$ correspond points $\widetilde{z_1}, \widetilde{z_2}$ of $\mO_{u_3}$. As $z_1 \sim z_2$ in $\chi'$, also $\widetilde{z_1} \sim \widetilde{z_2}$ in $\mS'$, a contradiction as $\widetilde{z_1}$ and $\widetilde{z_2}$ are points of an ovoid. Consequently $\chi$ does not contain triangles. 

From the foregoing it follows that $\chi$ is a GQ of order $(s,t)$ and that $\Gamma(\chi,\chi')^{\sigma^*} = \mE$. \eop \\

\section{Particular case}
\label{part}

Assume that $\gamma: \mA \mapsto \mE$ is a $\theta$-cover of $\mE$ (which is not necessarily the canonical projection). 
Let $x$ be a point of $\mS'$. The lines of $\mS$ not in $\mS'$ containing $x$ define $(t - t')/\theta$ rosettes by canonical projection, and these rosettes are covered by $t - t'$ mutually disjoint lines of $\mA$. The set of these $t - t'$ lines is a point $x^*$ of $\chi'$. 

\bl
The $t - t'$ lines of $\mA$ in $x^*$ are incident (in $\mS$) with a common point $x^{**}$ of $\mS'$.
\el 
{\em Proof}.\quad
Assume, by way of contradiction, that the distinct collinear points $x, y$ of $\mA$ define a common ovoid of $\mS'$ (through $\gamma$). This contradicts the fact that the $s$ points of the line $xy$ of $\mA$ define the $s$ distinct ovoids of a rosette in $\mS'$. 

Let $L, M$ be distinct lines of $\mA$ contained in $x^*$ and assume, by way of contradiction, that $L$ and $M$, considered as lines of $\mS$, are incident with distinct points $l$ and $m$ of $\mS'$. Let $y$ be a point of $L$, $z$ a point of $M$, with $y \sim z$. The ovoids of $\mS'$ defined by $y$ and $z$ belong to a common rosette $R$. Let $R(L)$ be the rosette defined by $L$ and let $R(M)$ be the rosette defined by $M$. 

First, assume that $R(L) = R(M)$. The ovoids defined by $y$ and $z$ belong to $R(L)$, $R(M)$ and $R$. Hence $R = R(L)$, a contradiction as $L$ and $yz$ are not disjoint in $\mA$. 

Next, assume that $R(L) \ne R(M)$. Let $O(y)$ be the ovoid defined by $y$ and let $O(z)$ be the ovoid defined by $z$. If $y'$ is any point subtending $O(y)$, then $y'$ is collinear in $\mS$ with at least one point $z'$ subtending $O(z)$. So in $\mS$ there arises a triangle with vertices $y', z'$ and $O(y) \cap O(z)$, a contradiction. 

We conclude that the $t - t'$ lines of $\mA$ in $x^*$ are incident with a common point $x^{**}$ of $\mS'$. 
\eop \\

\bt
\label{ident}
The generalized quadrangles $\chi'$ and $\mS'$ can be identified. 
\et
{\em Proof}.\quad
We may identify the point $x^*$ of $\chi'$ with the point $x^{**}$ of $\mS'$. Let $x^*$ and $y^*$ be distinct noncollinear points of $\chi'$. Assume by way of contradiction, that $x^{**}$ and $y^{**}$ are collinear in $\mS'$. Then the points $x = {(x^*)}^{\sigma^*}$ and $y = {(y^*)}^{\sigma^*}$ belong to some ovoid $\mO_z$ of $\mS'$. If $z'$ is a point of $\mC$ defining $\mO_z$, then $z' \sim x^{**}$ and $z' \sim y^{**}$ in $\mS$. hence there arises a triangle in $\mS$, a contradiction.
\eop \\

Theorem \ref{ident} leads to another proof of the lower decomposition property ($\sigma^*$ is the induced automorphism of $\mS'$). 

We keep using the same notation for the rest of this section.

\medskip
\subsection*{Condition (C)}

Consider again the geometry $\mE = \Gamma(\mS,\mS') = (\Omega,\Phi,\in)$. Let $M_1,M_2,\\
\ldots,M_r$, with $1 \leq r \leq \theta$, be lines  of $\mS$ covering a common rosette in $\mS'$ and let $L$ be a line of $\mS$ not in $\mS'$, which is not concurrent with $M_1,M_2,\ldots,M_r$ in $\mS$. Put $\mM := \{ M_1,M_2,\ldots,M_r \}$.
Assume that $x_0$ is the common point of $M_i$ and $\mS'$, with $i = 1,2,\ldots,r$, and that $x_1,x_2,\ldots,x_\alpha$ are  $\alpha$ points of $\mS \setminus \mS'$ on the lines $M_1,M_2,\ldots,M_r$, with $\alpha \in \{ s - 1, s\}$, and where for each $i$ the line $M_i$ contains at least one of the points $x_1,x_2,\ldots,x_\alpha$. Let $N_i$ be the line incident with $x_i$ and concurrent with $L$, with $i = 0,1,\ldots,\alpha$. Assume that all points $N_i \cap L$ are distinct, so that there arise $s$ or $s + 1$ points on $L$.
The set consisting of all point of $\mS'$ incident with the lines $N_0,N_1,\ldots,N_\alpha$ is denoted by $\mW_{L,\mM}$ if $N_0$ is not a line of $\mS'$ and by $\overline{\mW_{L,\mM}}$ if $N_0$ {\em is} a line of $\mS'$. We say that {\em Condition (C)} is satisfied if no sets $\mW_{L',\mM'}$ and $\mW_{L,\mM}$ with $\vert \mM' \vert = 1$ and $\vert \mM \vert > 1$, are isomorphic in $\mS'$, and if additionally no sets $\overline{\mW_{L',\mM'}}$ and $\overline{\mW_{L,\mM}}$, with $\vert \mM' \vert = 1$ and $\vert \mM \vert > 1$, are isomorphic in $\mS'$.\\

The next theorem shows that if $(C)$ holds, higher decomposition is possible.  

\bt
\label{thm8}
Consider a $\theta$-cover $\gamma: \mA \mapsto \mE$ and suppose that Condition (C) is satisfied. Then there exists an automorphism $\sigma$ of $\mS$ with $\sigma^*$ the restriction of $\sigma$ to $\mS'$.
\et
{\em Proof}.\quad
We will construct an automorphism $\sigma$ of $\mS$, such that $\sigma^*$ is the restriction of $\sigma$ to $\mS'$.

Choose a point $z$ of $\mA$ and choose a line $N$ of $\mS$ incident with $z$; the point of $\mS'$ incident with $N$ is denoted by $x^*$. The point $z$ covers an ovoid $\mO(z)$ of $\mS'$ and $\mO(z)$ is subtended by the points $z' = z_1,z_2,\ldots,z_{\theta}$ of $\mS$. Let $z^{\sigma} = z'$. Then the line $N^{\sigma}$ is the line $N' = z'x$, with $x = {x^*}^\sigma = {x^*}^{\sigma^*}$. 
For any line $T$ of $\mS'$, we define $T^\sigma = T^{\sigma^*}$. 

Now let $y$ be a point of $zx^*$, with $x^* \ne y \ne z$. The point $y$ covers an ovoid $\mO(y)$ of $\mS'$ and $\mO(y)$ is subtended by a unique point $y'$ of $N'$; let $y' = y^\sigma$. Next, consider a point $v$ of $\mA$, which is not collinear with $x^*$. Let $V$ be the line incident with $v$ and concurrent with $N$, 
let $w^*$ be the common point of $V$ and $\mS'$, and let $r$ be the common point of $N$ and $V$. The point $v$ covers and ovoid $\mO(v)$ of $\mS'$ and $\mO(v)$ is subtended by a point $v'$ of the line $V' = wr^\sigma$ (note that this is indeed a line), with ${w^*}^\sigma = w$. Let $V' = V^\sigma$ and $v' = v^\sigma$. Then $\sigma$ is defined for all points of $\mS$ not colllinear with $x^*$, for all points of the line $N$, for all points of $\mS'$, for the line $N$, for all lines concurrent with $N$ but not containing $x^*$, and for all lines of $\mS'$; for all these points and lines $\sigma$ preserves incidence.

Next, let  $U$ be a line of $\mS$ not concurrent with $N$ and not contained in $\mS'$. Let $u^*$ be the common point of $U$ and $\mS'$. The line $U$ of $\mA$ covers a rosette consisting of $s$ ovoids $\mO(u_1),\mO(u_2),\ldots,\mO(u_s)$ of $\mS'$, with $u_1,u_2,\ldots,u_s$ the points of $\mA$ incident with $U$. These ovoids have as common point $u = {u^*}^\sigma$. The ovoid $\mO(u_i)$ is subtended by the points $u_i^1,u_i^2,\ldots,u_i^{\theta}$ of the lines $U_1,U_2,\ldots,U_{\theta}$ incident with $u$, $i = 1,2,\ldots,s$. The points $u_i$ which are not collinear with $x^*$ are mapped by $\sigma$ onto points $u_i^{\sigma}$, where $u_i^{\sigma}$ is incident with a line $U_1,U_2,\ldots,U_{\theta}$. Call the subset of $\{U_1,U_2,\ldots,U_{\theta}\}$ which consists of the lines containing the points $u_i^{\sigma}$, with $u_i \not\sim x^*$ and $i \in \{1,2,\ldots,s\}$, $\mU'$, and put $r := \vert \mU'\vert$.
The lines concurrent with $N$ and $U$ which are not incident with $x^*$ are mapped by $\sigma$ onto the lines incident with $u_i^{\sigma}$, where $u_i \not\sim x^*$, and concurrent with $N^{\sigma}$, $i \in \{1,2,\ldots,s\}$, and the line incident with ${u^*}^{\sigma}$ and concurrent with $N^{\sigma}$. So there arise either $s$ or $s + 1$ lines intersecting $N^\sigma$ in distinct points.
These two sets of lines define isomorphic point sets $\mW_{N,\{U\}}$ and $\mW_{N^{\sigma},\mU'}$, or $\overline{\mW_{N,\{U\}}}$ and $\overline{\mW_{N^{\sigma},\mU'}}$ in $\mS'$. As Condition (C) is satisfied we have $r = 1$. So the points $u_i^{\sigma}$, with $u_i \not\sim x^*$, are incident with the line $U_j$ for some $j \in \{ 1,2,\ldots,\theta \}$, say $U_j = U_1$. Now we put $U^{\sigma} = U_1$.

We conclude that $\sigma$ is defined on all lines of $\mS$ not incident with $x^*$ and on all points not collinear with $x^*$. For all these points and lines $\sigma$ preserves incidence. One easily shows that such a $\sigma$ can be extended to an automorphism of $\mS$ for which $\sigma^*$ is the restriction to $\mS'$. 
\eop \\

\br{\rm
Note that the proof implies that there are $\theta$ such automorphisms $\sigma$. 
}\er

The first part of the next theorem yields the higher decomposition property for the case $\mS \cong \mQ(5,q)$, $\mS' \cong \mQ(4,q)$.

\bt
\begin{itemize}
\item[{\rm (i)}]
If $\mS = \mQ(5,q)$, $\mS' = \mQ(4,q)$, $\mC = \mA = \mS \setminus \mS'$, then there exists and automorphism $\sigma$ of $\mQ(5,q)$ with $\sigma^*$ the restriction of $\sigma$ to $\mQ(4,q)$. 
\item[{\rm (ii)}]
If $\mS = \mH(4,q^2)$, $\mS' = \mH(3,q^2)$, $\mC = \mA = \mS \setminus \mS'$, then there exists and automorphism $\sigma$ of $\mH(4,q^2)$ with $\sigma^*$ the restriction of $\sigma$ to $\mH(3,q^2)$. 
\end{itemize}
\et
{\em Proof}.\quad
By Theorem \ref{thm8} it is sufficient to show that Condition (C) is satisfied. In both cases sets $\mW_{L,\mM}$ and $\overline{\mW_{L,\mM}}$ are contained in a plane of the ambient projective space of $\mS'$, while sets $\mW_{L,\mM}$ and $\overline{W_{L,\mM}}$, with $\vert \mM \vert > 1$, are never contained in a plane of the ambient space of $\mS'$. It follows that Condition (C) is satisfied, and so the theorem is proved. 
\eop \\

\medskip
\section{Subtending in the Kantor-Knuth quadrangles}
\label{last}

Let $\mS$ be a nonclassical Kantor-Knuth GQ of order $(q,q^2)$. By \cite{Stab}, $\Aut(\mS)$ does not act transitively on the subquadrangles of order $q$ (which are all isomorphic to $\mQ(4,q)$). Let $\Gamma_1$ and $\Gamma_2$ be two such subGQs in different $\Aut(\mS)$-orbits (it is known that there are two such orbits).  If they are both doubly subtended, then we obtain two SPGs $\mE_1$ and $\mE_2$. Suppose they are isomorphic. 
Then if there would be an element $\kappa$ in $\Aut(\mS)$ sending $\mE_1$ to $\mE_2$, and so $\Gamma_1$ to $\Gamma_2$, we have a contradiction. 
So under this assumption (which we call (\#) for now), the SPGs are not isomorphic, and in \cite{Brown} Brown only found one particular subGQ which was doubly subtended. 
From this point of view it is of fundamental interest to check whether {\em all} $\mQ(4,q)$ subGQs of $\mS$ are doubly subtended. 

In this section we will show that the elements of one of these orbits do not have this property, so that no new SPGs can arise. 

Below, $\mS$ is a nonclassical Kantor-Knuth GQ of order $(q,q^2)$ | see \cite[\S 4.5 and \S 4.7.3]{TGQ}. It has a line $[\infty]$ which is fixed by $\Aut(\mS)$, and there are $q^3 + q^2$ subGQs of order $q$, all isomorphic to $\mQ(4,q)$, that contain the line $[\infty]$. There are no other subquadrangles of order $q$. All this can be found in \cite[\S 5.1]{TGQ}.

From \cite[Theorem 5.2.1(b)]{TGQ} the following result on orbits of subquadrangles can be deduced.

\medskip
\begin{itemize}
\item[{\bf Orbit property}]
There are two $\Aut(\mS)$-orbits of subquadrangles of order $q$; one of size $2q^2$ (called $\Omega_1$), one of size $(q - 1)q^2$ (called $\Omega_2$).
In particular, if $\Gamma$ is a line grid with parameters $(q,1)$ and containing $[\infty]$, then there are $2$ $\Omega_1$-subquadrangles containing $\Gamma$, and $q - 1$ $\Omega_2$-subquadrangles containing $\Gamma$. \\
\end{itemize}

In Brown \cite{Brown}, it is shown that the elements in $\Omega_1$ are all doubly subtended. Now suppose that the elements in $\Omega_2$ also are.
Consider any grid $\Gamma$ as above. As we assume that all $q + 1$ subquadrangles of order $q$ containing $\Gamma$ are doubly subtended, we obtain precisely $q + 1$ distinct involutory automorphisms of $\mS$, each fixing one of the subquadrangles of order $q$ containing $\Gamma$ pointwise. So the group $\vartheta$ generated by these involutions has at least size $q + 2$. Each element of $\vartheta$ fixes $\Gamma$ elementwise.
Now $\vartheta$ acts on the subquadrangles $\mS_1$ and $\mS_2$ in $\Omega_1$ that contain $\Gamma$. In particular, it follows that if $\theta_1$ is the involution that fixes $\mS_1$ elementwise, we have
\begin{equation}
2 \geq \vert \vartheta_{\mS_1}/\langle \theta_1 \rangle \vert \geq  \frac{q + 2}{4}, 
\end{equation}
since the subgroup of $\mQ(4,q)$ ($q$ odd) that fixes a grid with parameters $(q,1)$ elementwise, has size $2$. 

It follows that $q \leq 6$, so that $q$ is a prime. In those cases, $\mS \cong \mQ(5,q)$ \cite[\S 4.5]{TGQ}, contradiction. 

We have essentially proven the following.

\bt
No elements of $\Omega_2$ are doubly subtended. In particular, each element of $\Omega_2$ contains $(q + 1)q^2(q - 1)$ Kantor-Knuth ovoids which are $1$-subtended.
\et

{\em Proof}.\quad
The proof follows from the discussion preceding the statement of the theorem, plus the fact that for any subquadrangle of order $q$ containing $[\infty]$, any subtended ovoid is a Kantor-Knuth ovoid by \cite[Theorem 14]{Glasg}. \eop \\

In a forthcoming paper, we will study the assumption (\#), the Kantor-Knuth quadrangles and the associated SPGs in more depth.

\newpage
\appendix
\section{Extension in nonclassical cases}

In this appendix, we digress to focus on the extension property for $\mS'$ instead of $\mE$.

\subsection{TGQs}
A {\em translation generalized quadrangle} (TGQ) is a generalized quadrangle $\mS$ which has a point $u$ called {\em translation point} such that there is an abelian automorphism group $A \leq \Aut(\mS)_{[u]}$ that acts sharply transitively on the points not collinear with $u$. For details about the results and notions on TGQs used in this section, we refer to the monograph \cite{TGQ}.

Let $\mS$ be any thick TGQ of order $(q,q^2)$ which is not classical, and which has a subGQ $\mS'$ isomorphic to $\mQ(4,q)$ (if a GQ of order $(q,q^2)$ has a classical subGQ of order $(s',t')$ which is a geometrical hyperplane, then $s' = t' = q$, and the subGQ must be isomorphic to $\mQ(4,q)$). It is well known (see for example \cite{SFGQ}) that  $\mS'$ contains all the translation points of $\mS$ (and either there is one such point, or a line of translation points). Suppose the embedding 
$\iota: \mS' \hookrightarrow \mS$ has the higher extension property. Then any automorphism of the $\mQ(4,q)$-subquadrangle extends to an automorphism of $\mS$. In particular, it follows that $\Aut(\mS)_{\mS'}$ acts transitively on the points of $\mS'$, so any point of $\mS'$ is a translation point of $\mS$. So $\mS$ is classical, a contradiction. 

So any such $\mS$ gives rise to embeddings which do not have the higher extension property. Of the known examples, only the $\mQ(5,q)$ quadrangles and the Kantor-Knuth quadrangles also give rise to covers which are not $1$-covers (see, e.g., \cite{Stab}).

\subsection{Generalization to EGQs}

If we leave out the term ``abelian'' in the definition of TGQ, we have, by definition, an {\em elation generalized quadrangle} (EGQ) with {\em elation point} $u$. Most of the standard results on EGQs can be found in the book \cite{KTEGQ}.

Now consider the same situation as in the previous section, with ``TGQ'' replaced by ``EGQ.'' 
Suppose the embedding  $\iota: \mS' \hookrightarrow \mS$ has the higher extension property. 

Suppose first that $x$ is an elation point of $\mS$ which is contained in $\mS'$. 
Then as above any point of $\mS'$ is an elation point of $\mS$. So each point of $\mS$ is an elation point. By the main result of \cite{KTHVM}, $\mS$ is classical, a contradiction. 

Now suppose that $x$ is not contained in $\mS'$. Let $u \sim x$ be a point of $\mS'$. As $\mS'$ is classical, there is an automorphism $\kappa$ of $\mS'$ which moves $u$ to a collinear point $u' \ne u$. Let $\widetilde{\kappa}$ be an automorphism of $\mS$ which induces $\kappa$ in $\mS'$. Then $y := x^{\widetilde{\kappa}} \ne x$. If $y \sim x$, $yx$ is a line of elation points of $\mS$, so that $xy \cap \mS'$ is also an elation point of $\mS$, and then we are in the previous case. If $y \not\sim x$, each point of $\mS$ is an elation point, and then again by \cite{KTHVM}, $\mS$ is classical, contradiction.

\subsection{Quadrangles with reguli}

Let $\mS'$ be a GQ with finite parameters $(s,1)$, and let it be a full geometric hyperplane in a GQ $\mS$. Then $\mS$ necessarily has parameters $(s,s)$. Suppose that every automorphism of $\mS'$ extends to automorphisms of $\mS$. The automorphism group of $\mS'$ is isomorphic to $(\hS_{s + 1} \times \hS_{s + 1})\rtimes \mu_2$, where $\hS_\ell$ denotes the symmetric group on $\ell$ letters, and $\mu_2$ is a group of order $2$. Now if $s$ is a prime power, $\hS_{s + 1} \times \hS_{s + 1}$ (by which we mean the subgroup of $\Aut(\mS')$ that preserves the reguli of $\mS'$) 
contains the natural action(s) of $\PSL_2(s) \times \PSL_2(s)$ (in the language of De Kaey and Van Maldeghem \cite{DKVM}), and by {\em loc. cit.}, it follows that  $\mS \cong \mQ(4,s)$. However, in that case $\Aut(\mS)_{\mS'}$ induces $\mathbf{P\Gamma O}_4^{+}(s)$ (which has size $(s + 1)^2s^2(s - 1)^2\cdot 2h$, where $s = p^h$ with $p$ prime) on $\mS'$. It follows that $s \in \{ 2,3,4\}$. For $s = 2, 3$, one easily sees that the higher extension property is true. However, for $s = 4$, it is not: this follows immediately by comparing the sizes of $\mathbf{P\Gamma O}_4^{+}(4)$ and $(\hS_{5} \times \hS_{5})\rtimes \mu_2$, but can also easily be seen by observing that the action of the latter group would imply that every ovoid of $\mS'$ would be subtended by at least one point of $\mS$, and this is not possible in $\mQ(4,4)$. 

When $s$ is a not a prime power, see the discussion below, and Remark \ref{A.1}.\\

Suppose now that  $\Delta$ is an infinite generalized quadrangle which has a full thin subquadrangle $\Omega$ of order $(\omega,1)$ that is a geometrical hyperplane, and such that $\iota: \Omega \hookrightarrow \Delta$ has the higher extension property.

Let $\Delta$ and $\Omega$ be as above. Let $A$ be an automorphism group of $\Delta$ which induces the full automorphism group of $\Omega$ on $\Omega$. So $A$ induces $(\hS(\Omega_1) \times \hS(\Omega_2))\rtimes \mu_2$ on $\Omega$, where $\Omega_1$ and $\Omega_2$ are the two reguli of $\Omega$, and $\hS(S)$, with $S$ a set, denotes the symmetric group of $S$. Note that 
\begin{equation}
\vert A \vert \geq \vert \hS(\Omega_1) \vert = \vert \hS(\Omega_2) \vert > \vert \Omega_1 \vert = \vert \Omega_2 \vert = \omega.
\end{equation}

Now let $A_1$ be the subgroup of $A$ that fixes each line of $\Omega_2$; it induces $\hS(\Omega_1)$ on $\Omega_1$. 
Let $z$ be any point of $\Delta$ which is not contained in $\Omega$. For each line $M \in \Omega_2 \cup \Omega_1$, let $Z_M$ be the line incident with $z$ which meets $M$. Let $\mL$ be the set of lines $\{ Z_L \vert L \in \Omega_2\}$. Then as $\Omega$ is a geometrical hyperplane, we have that $t = \vert \mL \vert = \vert \Omega_2 \vert = \omega$, where $(\omega,t)$ is the order of $\Delta$. 

So the number of points of $\Delta$ equals $\omega$ as well. As $\vert A_1 \vert  > \omega$, we have $\vert {(A_1)}_z \vert > 1$. But the fixed elements structure of any element $\delta \ne \id$ of ${(A_1)}_z$ is a subquadrangle which is full and ideal, so by \cite[1.8.2]{POL}, $\delta$ is the identity in $A$, contradiction. 
So if $\Omega$ is supposed to be a geometrical hyperplane, then no examples can arise.

\br\label{A.1}
{\rm
Note that the latter approach also works for the finite case.}
\er

The next general question naturally comes to mind, but might be hard:

\begin{quote}
Classify the infinite quadrangles $\Delta$ of order $(\omega,t)$ with $t \geq \omega$ which have a full thin subquadrangle $\Omega$ of order $(\omega,1)$, such that $\iota: \Omega \hookrightarrow \Delta$ has the higher extension property. 
\end{quote}

\end{document}